\documentclass[12pt]{amsart}
\usepackage{xcolor}
\usepackage{tikz-cd}
\usepackage[margin=1in]{geometry}
\usepackage{todonotes} 
\usepackage{relsize}

\usepackage{amssymb, amsthm, array}
\usepackage{amsmath}

\usepackage[linktocpage=true, pdfencoding=auto, psdextra]{hyperref}

\newcommand{\ZK}{\mathcal{Z}_K} %moment-angle complex
\newcommand{\KJ}{K_J} %full subcomplex of K on J
\newcommand{\A}{\mathcal{A}}
\newcommand{\F}{\mathbb{F}_2}
\renewcommand{\star}{\mathrm{star}}
\newcommand{\link}{\mathrm{link}}
\newcommand{\Cone}{\mathrm{Cone}\,}

\newtheorem{thm}{Theorem}[subsection]

\theoremstyle{plain}
\newtheorem{prop}[thm]{Proposition}
\newtheorem{lem}[thm]{Lemma}
\newtheorem{cor}[thm]{Corollary}

\theoremstyle{definition}
\newtheorem{defn}[thm]{Definition}

\theoremstyle{remark}
\newtheorem{rem}[thm]{Remark}
\newtheorem{eg}[thm]{Example}

\unitlength=1mm \linethickness{0.5pt} 

\def\F{\mathbb F}

\def\R{\mathbb R}

\def\cZ{\mathcal Z}

\begin{document}

\title[Steenrod operations on polyhedral products]{Steenrod operations on polyhedral products}

\author[Sanjana Agarwal et al.]{Sanjana Agarwal}
\address[Sanjana Agarwal]{Department of Mathematics\\ Indiana University\\ Bloomington, IN 474005\\ USA}
\email{sanjagar@iu.edu}

\author[]{Jelena Grbi\' c} 
\address[Jelena Grbi\' c]{School of Mathematical Sciences, University of Southampton, Southampton, SO17  1BJ, UK}
\email{j.grbic@soton.ac.uk}

\author[]{Michele Intermont}
\address[Michele Intermont]{Department of Mathematics\\ Kalamazoo College\\ Kalamazoo, MI 49006\\ USA}
\email{intermon@kzoo.edu}

\author[]{Milica Jovanovi\' c}
\address[Milica Jovanovi\' c]{University of Belgrade\\ Faculty of mathematics\\ Studentski trg 16\\ Belgrade\\ Serbia}
\email{milica.jovanovic@matf.bg.ac.rs}

\author[]{Evgeniya Lagoda} 
\address[Evgeniya Lagoda]{Freie Universit\"at Berlin \\ Institute of Mathematics \\ Arnimallee 2\\ 14195 Berlin, Germany }
\email{evgeniya.lagoda@fu-berlin.de}

\author[]{Sarah Whitehouse}
\address[Sarah Whitehouse]{School of Mathematics and Statistics\\
University of Sheffield\\ Sheffield\\ S3 7RH\\ UK}
\email{s.whitehouse@sheffield.ac.uk}

\date{\today}



\begin{abstract}
We describe the action of the mod $2$ Steenrod algebra on the cohomology of various polyhedral products and related spaces. We carry this out for Davis-Januszkiewicz spaces and their generalizations, for moment-angle complexes as well as for certain polyhedral joins. By studying the combinatorics of underlying simplicial complexes, we deduce some consequences for the lowest cohomological dimension in which non-trivial Steenrod operations can appear. 

We present a version of cochain-level formulas for Steenrod operations on simplicial complexes. We explain the idea of ``propagating'' such formulas from a simplicial complex $K$ to polyhedral joins over $K$ and we give examples of this process. We tie the propagation of the Steenrod algebra actions on polyhedral joins to those on moment-angle complexes. Although these are cases where one can understand the Steenrod action via a stable homotopy decomposition, we anticipate applying this method to cases where there is no such decomposition.
\end{abstract}

\maketitle

\tableofcontents

\setcounter{section}{0} %%%%%%%%%%%%%%%%%%%%%%%%%%%%%%%%%%%%%%%%%%%%%%%%%%%%%%%%%%%
\section{Introduction}
%%%%%%%%%%%%%%%%%%%%%%%%%%%%%%%%%%%%%%%%%%%%%%%%%%%%%%%%%%%

Beginning with a simplicial complex $K$ on $m$ vertices and $m$ topological pairs $(\underline X, \underline A)=((X_1, A_1),\ldots, (X_m, A_m))$, the \textit{polyhedral product} $(\underline X, \underline A)^K$ is a colimit construction over the face category of $K$, that is,
\[
(\underline X, \underline A)^K=\bigcup_{\sigma\in K} (\underline X, \underline A)^\sigma\subset X_1\times\ldots\times X_m
\]
where $(\underline X, \underline A)^\sigma=\{(x_1,\ldots, x_m)\in \prod X_i \ | \ x_i\in A_i \text{ if } i\notin \sigma\}$.

This construction covers a surprisingly large family of topological spaces. For example, when $(\underline X, \underline A)=(D^n, S^{n-1})$ and $K$ is a triangulation of a sphere, the polyhedral product $(D^n, S^{n-1})^K$ is a manifold. In particular, the manifold $(D^2, S^1)^K$, called a moment-angle manifold, comes equipped with a $T^m$-action and is a topological analogue of projective toric varieties studied in algebraic geometry.  It is also an intersection of quadrics which is an object of interest in complex geometry. For an arbitrary simplicial complex $K$, the  moment-angle complex $(D^2, S^1)^K$ can be identified as the complement of a complex coordinate subspace arrangement, linking polyhedral products to combinatorics and commutative algebra. 

The main goal of this paper is to introduce a systematic study of the action of the Steenrod algebra on the cohomology of polyhedral products. So far sporadic observations regarding the Steenrod algebra actions in toric topology have been made in several papers (see for example~\cite{MR1608269},~\cite{MR3120631}). Here we describe the action of the mod $2$ Steenrod algebra~$\mathcal{A}$ on the cohomology of polyhedral products and related spaces. For the case of Davis-Januszkiewicz spaces
there is a well-known explicit presentation of the cohomology
as the Stanley-Reisner algebra and this has an action of the Steenrod algebra, stemming from the action on complex projective space $\mathbb {C}P^{\infty}$. 
We also note a generalization to polyhedral products $(\underline{X}, \underline {*})^K$ that works in a similar way.  

For moment-angle complexes $(D^2, S^1)^K$, the situation is more complicated.
This is because, while the cohomology can be described as  a $\rm{Tor}$ algebra, there is no known presentation of the cohomology in terms of algebra generators and relations. However, the Steenrod action can be understood via the Hochster splitting and the underlying homotopy decomposition of moment-angle complexes after suspension. As an example we consider for $K$ a triangulation of real projective space $\mathbb{R}P^2$ and 
we deduce some consequences for the lowest dimension in which a non-trivial $Sq^1$ can appear in the cohomology of a moment-angle complex. When $K$ is a graph, a 1-dimensional simplicial complex, we deduce that there is no non-trivial Steenrod operation on the moment-angle complex $(D^2, S^1)^K$.

 Consideration of the Hochster splitting motivates us to understand the action of the Steenrod algebra on simplicial complexes. So we also present a combinatorial formula for the action of Steenrod squares on the cochains of the simplical complex $K$ which is a version of formulas found in~\cite{GONZALEZDIAZ199989}. Then we examine how the action of the Steenrod algebra on the cohomology of $K$ can be propagated to more complicated simplicial complexes built out of $K$. That is, given $K$ and non-trivial actions of Steenrod operations on $K$, 
 we study operations on various polyhedral join constructions built from $K$. Finally, we relate the results back to our earlier analysis by considering the moment-angle complexes associated to these complexes.
\medskip

The paper is organized as follows. In section 2 we recall the necessary definitions and notation. Section 3 describes the action of the Steenrod algebra $\A$ on the cohomology of Davis-Januszkiewicz spaces $(\mathbb{C}P^{\infty},*)^K$ and on the more general polyhedral products of the form $(\underline{X}, \underline{*})^K$. In section 4 we turn our attention to describing the action of $\mathcal{A}$ on the cohomology of moment-angle complexes. We also provide some results when the simplicial complex $K$ contains the minimal triangulation of $\mathbb{R}P^2$.
Section 5 covers the combinatorial formula for computing the action of $Sq^n$ on cochains. In Section 6 we turn our attention to polyhedral joins, in particular to substitution complexes $K\langle K_1, \dots , K_m\rangle$ and obtain some results on how Steenrod operations on $K, K_1, \dots , K_m$ may be used to understand
Steenrod operations on the substitution complex itself. Finally we use this to deduce some information about how Steenrod operations act on the associated moment-angle complex.

\subsection*{Acknowledgements} We would like to thank the
Hausdorff Research Institute for Mathematics for hosting the Women in Topology IV workshop and for financial support. We would also like to thank the
Foundation Compositio Mathematica, the Foundation Nagoya Mathematical Journal and the K-theory Foundation for financial support for this event.

\setcounter{section}{1} %%%%%%%%%%%%%%%%%%%%%%%%%%%%%%%%%%%%%%%%%%%%%%%%%%%%%%%%%%%
\section{Preliminaries}
%%%%%%%%%%%%%%%%%%%%%%%%%%%%%%%%%%%%%%%%%%%%%%%%%%%%%%%%%%%

In this section we review some necessary background.  We begin with the construction of the polyhedral product, before giving some brief
background on Steenrod operations.

\subsection{Notation} We denote by
${[ m ]}$ the set $\{1,\dots,m\}$ endowed with its natural order and we use this to label the vertices of a simplicial complex, $K$, on $m$ vertices.
 The full subcomplex of $K$ obtained by restricting to $J \subseteq [m]$ will be denoted by $K_{J}$.

We work with cohomology with $\F_2$ coefficients throughout the paper and we will write
$H^*(-)$ for $H^*(-; \mathbb{F}_2)$. Where we discuss other coefficients, we will write them explicitly. We write $\A$ for the mod $2$ Steenrod algebra.

Suspension for both spaces and graded modules is denoted by $\Sigma$. For suspension of a graded module $M$, we have $(\Sigma M)^n=M^{n-1}$. 

\subsection{Polyhedral products}
\begin{defn}
Given a simplicial complex $K$ on $m$ vertices, and $m$ topological pairs $(\underline X, \underline A)=((X_1, A_1),\ldots, (X_m, A_m))$, the \emph{polyhedral product} $(\underline X, \underline A)^K$ is the subspace of the product
$\prod_{i=1}^m X_i$ given by
\[
(\underline X, \underline A)^K=\bigcup_{\sigma\in K} (\underline X, \underline A)^\sigma\]
where \[(\underline X, \underline A)^\sigma=\prod_{i=1}^m Y_i\quad {\rm \ and \ }\quad  Y_i = \begin{cases} X_i &{\rm \ if \ } i \in \sigma \\ A_i & {\rm \ if \ } i \notin \sigma.\end{cases}\]
\end{defn}

There are some special cases of this construction worth mentioning.
When $K$ is an $(m-1)$-simplex, then $(\underline X, \underline A)^K = X_1 \times \dots \times X_m$.
When $K$ is the disjoint union of $m$ points and the spaces $A_i$ are all taken as a point, then $(\underline X, \underline *)^K$ is the wedge $X_1 \vee \dots \vee X_m$. We write $(X,A)^K$ when $(\underline X, \underline A)=((X_1, A_1),\ldots, (X_m, A_m)) =((X, A),\ldots, (X, A))$. The Davis-Januszkiewicz space, $DJ(K) = (\mathbb{C}P^{\infty}, *)^K$, and the moment-angle complex, $\ZK = (D^2, S^1)^K$, are key examples in this paper.

\subsection{Steenrod operations}
Steenrod operations are the primary stable operations of ordinary cohomology with $\mathbb{F}_p$ coefficients for a prime $p$. In this paper, we will focus on the case $p=2$, but many of the results have obvious analogues for odd primes.
In particular this applies to results which come from stable homotopy decompositions.

For a space $X$, its cohomology with $\mathbb{F}_2$ coefficients, which we denote $H^*(X)$, has an action of the mod $2$ Steenrod algebra $\mathcal{A}$. 
This action satisfies a stability property, meaning that
it commutes with the suspension isomorphisms.
The $\mathbb{F}_2$-algebra $\mathcal{A}$ has generators $Sq^n$ for $n\geq 1$ subject to the Adem relations. The interaction of the $\mathcal{A}$-algebra structure and the cup product is controlled by the Cartan formula.
This can be conveniently captured by the statement that the total Steenrod operation,
that is, the formal sum $\sum_{n=0}^\infty Sq^n$, is a ring map. Here $Sq^0=1$, acting as the identity operation. The instability condition satisfied by the cohomology of a space is that if $x$ is a class in degree $d$ then 
$Sq^d(x)=x^2$ and $Sq^n(x)=0$ for $n>d$.
This ensures that the total Steenrod operation acts as a finite sum on each homogeneous element.

An introductory account of mod $2$ Steenrod operations and their properties can be found in~\cite[Section 4L]{Hatcher}. Classical treatments include ~\cite{MosherTangora} and~\cite{Steenrod} .

\setcounter{section}{2} %%%%%%%%%%%%%%%%%%%%%%%%%%%%%%%%%%%%%%%%%%%%%%%%%%%%%%%%%%%
\section{Davis-Januszkiewicz spaces and generalizations}
%%%%%%%%%%%%%%%%%%%%%%%%%%%%%%%%%%%%%%%%%%%%%%%%%%%%%%%%%%%

In this section we describe the action of the Steenrod algebra $\mathcal{A}$ on the cohomology of
Davis-Januszkiewicz spaces via the description of the cohomology as the
Stanley-Reisner algebra. We also generalize this result to polyhedral products $(\underline{X}, \underline{*})^K$ for arbitrary pointed topological spaces $X_i$.

\subsection{Davis-Januszkiewicz spaces}

To a simplicial complex $K$ on $[m]$ and a commutative ring $R$, one can associate the Stanley-Reisner algebra 
\begin{equation}
    \label{SR}
R[K] =\frac{R[x_1, x_2, \dots, x_m]}{I_K}
\end{equation}
where the Stanley-Reisner ideal $I_K$ is generated by the squarefree monomials corresponding to non-faces of $K$,
\[
    I_K=\langle x_{i_1}x_{i_2}\dots x_{i_k}\,|\,(i_1, i_2, \dots, i_k)\notin K \rangle.
\]
Recall~\cite[Proposition 4.3.1]{buchstaber2014toric} that the cohomology algebra of the Davis-Januszkiewicz space $DJ_K$ is the Stanley-Reisner algebra
\begin{equation}
\label{cohmlgyDJK}
	H^*(DJ_K; R)\cong R[K]	
\end{equation}
where each $x_i$ is of degree 2.

Since our goal is to study the action of the mod $2$ Steenrod algebra $\mathcal{A}$,
we now specialize to the case $R=\mathbb{F}_2$.

We define an $\mathcal{A}$-algebra structure on the polynomial algebra $\mathbb{F}_2[x_1, x_2, \dots, x_m]$ by requiring the total Steenrod operation $Sq$ to act on generators by $Sq(x_i)=x_i+x_i^2$, and extending to $\mathbb{F}_2[x_1, x_2, \dots, x_m]$ by requiring $Sq$ to be an $\mathbb{F}_2$-algebra map. Of course, as discussed below, this is the
usual action on the cohomology of a product of infinite complex projective spaces.

\begin{lem}\label{Case1Lemma2}
Let $K$ be a simplicial complex on $[m]$. An $\mathcal A$-algebra structure on the Stanley-Reisner algebra $\mathbb{F}_2[K]$ can be defined on generators by $Sq(x_i)=x_i+x_i^2$.
\end{lem}

\begin{proof} 
Given the definition of the $\mathcal A$-algebra structure on the polynomial algebra, we only need to check that the Stanley-Reisner ideal $I_K$ is an $\mathcal{A}$-ideal.
For an arbitrary squarefree monomial of the form $x_{i_1}\cdots x_{i_k}$ for which $(i_1,\ldots,i_k)\notin K$, we show that $Sq(x_{i_1}\cdots x_{i_k})$ belongs to $I_{K}$. This follows by expanding the product
$Sq(x_{i_1}\cdots x_{i_k})=\prod_{j=1}^k Sq(x_{i_j})=\prod_{j=1}^k (x_{i_j}+x_{i_j}^2)$ into a sum of monomials in which every summand is  divisible by $x_{i_1}\cdots x_{i_k}$.
\end{proof}

To understand the action of the Steenrod algebra on $H^*(DJ_K)$, we consider isomorphism~\eqref{cohmlgyDJK} in more detail. The complex projective space $\mathbb CP^\infty$ can be given a CW structure with exactly one cell in each even dimension. We endow the product  $(\mathbb CP^\infty)^m$ with the product CW structure which consequently consists of only even-dimensional cells. Since $DJ_K$ is a cell subcomplex of $(\mathbb CP^\infty)^m$, the same holds for $DJ_K$. This implies that the cohomology algebras of $(\mathbb CP^\infty)^m$ and $DJ_K$ coincide with their cellular cochain algebras.

The inclusion $i\colon DJ_K\hookrightarrow (\mathbb CP^\infty)^m$ induces a homomorphism of differential graded algebras on cellular cochains
\[\mathcal C^*((\mathbb CP^\infty)^m)\longrightarrow \mathcal C^*(DJ_K).\]
As an $\mathbb F_2 $-vector space $\mathcal C^*((\mathbb CP^\infty)^m)$ has a basis of cochains $(D_{j_1}^{2k_1}\cdots D_{j_p}^{2k_p})^*$ dual to the products of cells $D_{j_1}^{2k_1}\times \cdots \times D_{j_p}^{2k_p}$.
Note that the cup product of basis elements is again a basis element.

Denote by $\Phi$ the $\mathcal A$-algebra isomorphism $\mathcal C^*((\mathbb CP^\infty)^m)\longrightarrow \mathbb F_2[x_1,x_2,\ldots,x_m]$.
The isomorphism $\Phi$ maps the cochain $(D_{j_1}^{2k_1} \cdots D_{j_p}^{2k_p})^*$ to the monomial $x_{j_1}^{k_1}\cdots x_{j_p}^{k_p}$.

\begin{lem}
The isomorphism~\eqref{cohmlgyDJK} is an $\mathcal A$-algebra isomorphism.
\end{lem}

\begin{proof}

We denote the algebra isomorphism $\eqref{cohmlgyDJK}$ by $\Psi$. Via $\Psi$, the generator  $(D_{j_1}^{2k_1}\cdots D_{j_p}^{2k_p})^*$ of $H^*(DJ_K)$ is identified with the equivalence class of the monomial $x_{j_1}^{k_1}\cdots x_{j_p}^{k_p}$ in the Stanley-Reisner algebra $\mathbb F_2[x_1,x_2,\ldots,x_m]/{I_K}$. 
Considering the $\mathcal A$-algebra structure on $\mathbb F_2[K]$ given in Lemma~\ref{Case1Lemma2}, we show that $\Psi$ is an $\mathcal{A}$-algebra isomorphism.

By the surjectivity of $i^*\colon H^*((\mathbb CP^\infty)^m)\longrightarrow H^*(DJ_K)$, by the naturality of Steenrod operations, and since $\Phi$ is an $\mathcal{A}$-algebra isomorphism, we have 
\begin{align*}
\Psi\circ Sq(D_{j_1}^{2k_1}\cdots D_{j_p}^{2k_p})^*&=\Psi\circ i^*\circ Sq(D_{j_1}^{2k_1}\cdots D_{j_p}^{2k_p})^*= [\Phi\circ Sq(D_{j_1}^{2k_1}\cdots D_{j_p}^{2k_p})^*]\\
&=[Sq\circ \Phi(D_{j_1}^{2k_1}\cdots D_{j_p}^{2k_p})^*]=[Sq(x_{j_1}^{k_1}\cdots x_{j_p}^{k_p})].
\end{align*}
At the same time, by the definition of the Steenrod action on $\mathbb F_2[x_1,x_2,\ldots,x_m]/{I_K}$,
\[
Sq([x_{j_1}^{k_1}\cdots x_{j_p}^{k_p}])=[Sq(x_{j_1}^{k_1}\cdots x_{j_p}^{k_p})] .
\]
Thus,
\[
\Psi\circ Sq(D_{j_1}^{2k_1}\cdots D_{j_p}^{2k_p})^*=Sq\circ \Psi(D_{j_1}^{2k_1}\cdots D_{j_p}^{2k_p})^*
\]
which proves that $\Psi$ is an $\mathcal A$-algebra isomorphism.
\end{proof}

\subsection{Generalizations: $(\underline{X}, \underline{*})^K$}

Davis-Januszkiewicz spaces can be generalized to the polyhedral products $(\underline{X}, \underline{*})^K$ for arbitrary pointed topological spaces $X_i$.  Recall~\cite[Theorem 2.35]{MR2673742} that if $R$ is a ring such that the natural map 
\[
H^*(X_{j_1}; R)\otimes\cdots \otimes H^*(X_{j_k}; R)\longrightarrow H^*(X_{j_1}\times\cdots\times X_{j_k}; R)
\]
is an isomorphism for any $\{ j_1,\ldots, j_k\}\subset [m]$, then there is an isomorphism of algebras
\begin{equation}
    \label{chlgyX*K}
H^*((\underline {X}, \underline *)^K; R)\cong\frac{\bigotimes_{i=1}^m H^*(X_i; R)}{I_K(X_1,\dots, X_m)}
\end{equation}
where $I_K(X_1,\dots, X_m)=\langle  x_{i_{1}}x_{i_2}\dots x_{i_k}\,|\, x_{i_j}\in \widetilde{H}^*(X_{i_j}; R) \text{\ and\ }(i_1, i_2, \dots, i_k)\notin K \rangle$.

Specialising to the case $R=\mathbb{F}_2$, our goal is to show that isomorphism~(\ref{chlgyX*K}) is an $\mathcal A$-algebra isomorphism. The proofs follow the same pattern as for the case of Davis-Januszkiewicz spaces.
The action of the Steenrod algebra $\mathcal A$ on $\otimes_{i=1}^m H^*(X_i)/I_K(X_1,\dots, X_m)$ is induced from that on $\otimes_{i=1}^m H^*(X_i)$ which as a tensor product of $\mathcal A$-modules  has a standard $\mathcal A$-module structure coming from the comultiplication on $\mathcal A$.
\begin{lem}\label{Case2Lemma2}
The ideal $I_K(X_1,\dots, X_m)$ is an $\mathcal{A}$-ideal.
\end{lem}
\begin{proof}
For $x_{i_1}x_{i_2}\cdots x_{i_k}\in I_K(X_1, \dots , X_m)$, that is, for $(i_1,i_2,\ldots, i_k)\notin K$ and $x_{i_j}\in \widetilde{H}^*(X_{i_j};\mathbb F_2)$ for each $j\in\{1,\ldots,k\}$, we prove that $Sq (x_{i_1}x_{i_2}\cdots x_{i_k}) \in I_K(X_1,\dots, X_m)$. As the total Steenrod operation is multiplicative, 
\[
 Sq(x_{i_1}x_{i_2}\cdots x_{i_k}) = Sq(x_{i_1})Sq(x_{i_2})\cdots Sq(x_{i_k}). 
\]
Since $Sq(x_{i_j})\in \widetilde{H}^*(X_{i_j})$, it follows that $ Sq(x_{i_1}x_{i_2}\cdots x_{i_k})$ is a sum of monomials of the form $y_{i_1}y_{i_2}\cdots y_{i_k}$, where each $y_{i_j}\in \widetilde{H}^*(X_{i_j})$. Since $(i_1,\ldots,i_k)\notin K$, we obtain that $y_{i_1}y_{i_2}\cdots y_{i_k}\in I_K(X_1,\dots, X_m)$ and thus $Sq(x_{i_1}x_{i_2}\cdots x_{i_k})\in I_K(X_1,\dots, X_m)$. 
\end{proof}

\begin{lem}
The isomorphism~\eqref{chlgyX*K} is an $\mathcal A$-algebra isomorphism.
\end{lem}
\begin{proof}
Denote isomorphism~\eqref{chlgyX*K} by $\Psi$. By Lemma \ref{Case2Lemma2}, since $I_K(X_1,\dots, X_m)$ is an $\mathcal{A}$-ideal, the $\mathcal{A}$-algebra structure on the tensor product $\otimes_{i=1}^m H^*(X_i)$ passes to the quotient ring $\otimes_{i=1}^m H^*(X_i)/{I_K(X_1,\dots, X_m)}$. We prove that $\Psi$ is also an $\mathcal{A}$-algebra isomorphism.

Consider the following commutative diagram
\begin{equation}\label{eq:A-algebra diagram}
\begin{tikzcd}
H^*(X_1\times \ldots \times X_m) \arrow[d, "i^*"] \arrow[rr, "\Phi"]       &  & \underset{i=1}{\overset{m}{\bigotimes}} H^*(X_i) \arrow[d]                         \\
{H^*((\underline{X},\underline{*})^K)} \arrow[rr, "\Psi"] &  & {\underset{i=1}{\overset{m}{\bigotimes}} H^*(X_i)/{I_K(X_1,\dots, X_m)}} 
\end{tikzcd}  
\end{equation}
where $i\colon (\underline{X},\underline{*})^K\to X_1\times \ldots \times X_m$ is the inclusion and the right hand arrow is the projection to the quotient.
Let $[\alpha]\in H^*((\underline{X},\underline{*})^K)$ be a class that is identified under the isomorphism $\Psi$ with $[x_{j_1}\cdots x_{j_p}]$ in  $\otimes_{i=1}^m H^*(X_i)/{I_K(X_1,\dots, X_m)}$. By the commutativity of~\eqref{eq:A-algebra diagram}, there is a class $[\Tilde{\alpha}] \in (i^*)^{-1}([\alpha])\subseteq H^*(X_1\times \ldots \times X_m)$ that is mapped under the Künneth isomorphism $\Phi$ to $x_{j_1}\cdots x_{j_p} \in \otimes_{i=1}^m H^*(X_i)$. 
Then 
\[
\Psi\circ Sq([\alpha])=\Psi\circ i^*\circ Sq([\Tilde{\alpha}])=[\Phi\circ Sq([\Tilde{\alpha}])]=[ Sq\circ\Phi([\Tilde{\alpha}])]= [Sq(x_{j_1}\cdots x_{j_p})]
\]
where the first equality is by naturality of Steenrod operations and the second equality is by commutativity of~\eqref{eq:A-algebra diagram}.
At the same time, by the definition of the Steenrod action on $\otimes_{i=1}^m H^*(X_i)/{I_K(X_1,\dots, X_m)}$,
\[
Sq([x_{j_1}\cdots x_{j_p}])=[Sq(x_{j_1}\cdots x_{j_p})] .
\]
Thus,
\[
\Psi\circ Sq([\alpha])=Sq\circ \Psi([\alpha]).\qedhere
\]
\end{proof}

\begin{rem}
    The same proof can be extended to the cohomology splittings of~\cite[Theorem~3.3, Theorem 3.4]{MR3235796}, showing that the isomorphisms hold as algebras over the Steenrod algebra.
\end{rem}

\begin{eg}
For $(X,*)= (\mathbb RP^{\infty}, *)$ consider the polyhedral product $(\mathbb RP^{\infty}, *)^{K}$. 
Recall that $H^*(\mathbb RP^{\infty})\cong \mathbb F_2[x]$, where $\deg(x)=1$, and the action of the total Steenrod square on the generator is given by $Sq(x)=x+x^2$. 
Thus, by~\eqref{chlgyX*K}, the mod-2 cohomology of $(\mathbb RP^{\infty}, *)^{K}$ is isomorphic to the Stanley–Reisner ring  $\mathbb F_2[K]=\mathbb F_2[x_1, \ldots, x_m]/I_K$. 
It follows that the action of the Steenrod squares on $H^*((\mathbb RP^{\infty}, *)^{K})$ is described by $Sq([x_i])=[Sq(x_i)]=[x_i+x_i^2]$.
\end{eg}

\setcounter{section}{3} %%%%%%%%%%%%%%%%%%%%%%%%%%%%%%%%%%%%%%%%%%%%%%%%%%%%%%%%%%%
\section{Moment-angle complexes}
%%%%%%%%%%%%%%%%%%%%%%%%%%%%%%%%%%%%%%%%%%%%%%%%%%%%%%%%%%%
Describing the Steenrod operations for moment-angle complexes $\ZK$ is more complicated than for Davis-Januszkiewicz spaces because there is no known presentation of the cohomology of $\ZK$. Instead, for $K$ on $[m]$ the cohomology as an algebra is given by
\[
H^*(\mathcal{Z}_K; \mathbb{F}_2)
	\cong {\rm{Tor}}_{\mathbb{F}_2[x_1, x_2, \dots, x_m]}(\mathbb{F}_2[K],\mathbb{F}_2)
\]
see~\cite[Theorem 4.5.4]{buchstaber2014toric}.
At the moment, we do not know how to describe the Steenrod action directly on the Tor; we leave this for a future project.
There is, however, a stable splitting of $\ZK$ which we use to obtain some understanding of the Steenrod operations.

\subsection{Splitting after suspension}

In one sense, complete information about the Steenrod operations on moment-angle complexes can be read off from the homotopy decomposition 
\begin{equation}
    \label{decomp}
\Sigma\ZK\simeq \bigvee_{J\notin K} \Sigma^{|J|+2} |\KJ|
\end{equation}
 see~\cite[Corollary 8.3.6(b)]{buchstaber2014toric}. Note that the wedge here could be expressed in terms of all $J\subseteq [m]$, or as we have done, just such $J\notin K$, since $|\KJ|$ is contractible if $J\in K$.

\begin{prop}
\label{prop:splitting}
There is an isomorphism of $\mathcal{A}$-modules
\begin{equation}
    \label{cohlgy_splitting}
\widetilde{H}^*(\ZK)
\cong\bigoplus_{J\notin K} \Sigma^{|J|+1}\widetilde{H}^{*}(\KJ).
\end{equation}
\end{prop}

\begin{proof}
From the homotopy decomposition~\eqref{decomp}, there
is an additive splitting of cohomology
\[
\widetilde{H}^*(\ZK)
\cong\bigoplus_{J\notin K} \widetilde{H}^{*}( \Sigma^{|J|+1}\KJ).
\]
Since the Steenrod operations are stable, this gives the stated isomorphism of $\A$-modules.
\end{proof}

Thus the only non-trivial Steenrod operations are internal to the summands in this cohomological decomposition. This
is in contrast to the cup product structure: cup products within summands are of course trivial since each summand is the cohomology of a suspension; the only potentially non-trivial cup products pair the summands corresponding to full subcomplexes on $J$ and $L$
to the summand corresponding to $J\cup L$ if $J\cap L=\emptyset$, see~\cite[Proposition 3.2.10]{buchstaber2014toric}.

\subsection{Main example}\label{sec:an_example}
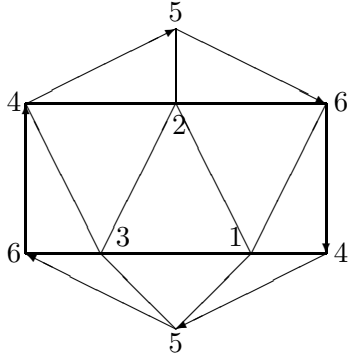
\begin{figure}[ht]
  \begin{center}
  \begin{picture}(120,40)
  \put(59,-3){\small 5}
  \put(59,41){\small 5}
  \put(81,29){\small 6}
  \put(81,9){\small 4}
  \put(37.5,29){\small 4}
  \put(37.5,9){\small 6}
  \put(59.5,26){\small 2}
  \put(52,11){\small 3}
  \put(67,11){\small 1}
  \put(60,0){\vector(-2,1){20}}
  \put(40,10){\vector(0,1){20}}
  \put(40,30){\vector(2,1){20}}
  \put(60,40){\vector(2,-1){20}}
  \put(80,30){\vector(0,-1){20}}
  \put(80,10){\vector(-2,-1){20}}
  \put(40,10){\line(1,0){40}}
  \put(40,30){\line(1,0){40}}
  \put(60,0){\line(-1,1){10}}
  \put(60,0){\line(1,1){10}}
  \put(50,10){\line(-1,2){10}}
  \put(50,10){\line(1,2){10}}
  \put(70,10){\line(-1,2){10}}
  \put(70,10){\line(1,2){10}}
  \put(60,30){\line(0,1){10}}
  \end{picture}
  \end{center}
  \caption{$P^2_6$, a $6$-vertex triangulation of $\mathbb R P^2$. Here the vertices with the same label are identified, and the
boundary edges are identified according to the orientation shown.}
  \label{rptri}
\end{figure}

Let us consider a specific example which combinatorially detects $Sq^1$ on moment-angle complexes. For the topological space $\mathbb{R}P^2$, we know
 $Sq^1 \colon H^1(\mathbb{R}P^2) \to H^2(\mathbb{R}P^2)$ is non-trivial.  Let $K$ be the triangulation of $\mathbb{R}P^2$ on six vertices, $P^2_6$, given by Figure~\ref{rptri}. Given that $P^2_6$ has a non-trivial $Sq^1$, we describe the action of $Sq^1$ on the moment-angle complex over $P^2_6$.

By the homotopy splitting \eqref{decomp}, 
\begin{align*}
  \label{stablerp2example}
  \Sigma\cZ_{P^2_6} &\simeq (\Sigma^{5}S^1)^{\vee 10}\vee(\Sigma^{6}S^1)^{\vee 15}\vee(\Sigma^{7}S^1)^{\vee 6}
  \vee\Sigma^8\R P^2\\
   &\simeq (S^6)^{\vee 10} \vee (S^7)^{\vee 15} \vee (S^8)^{\vee 6} \vee \Sigma^8\R P^2
\end{align*}
where $X^{\vee k}$ denotes the $k$-fold wedge of~$X$. Here we are using the fact that all the strict subcomplexes of $P^2_6$ are homotopic to the circle, and $P^2_6$ itself is homotopic to $\R P^2$.
For this example we have a stronger result,  namely an unstable splitting~\cite{MR3461047}
\begin{equation}
\label{Unstable_splitting}
  \cZ_{P^2_6} \simeq (S^5)^{\vee 10}\vee(S^6)^{\vee 15}\vee(S^7)^{\vee 6}
  \vee\varSigma^7\R P^2.
\end{equation}

Thus, the cohomology of $\cZ_{P^2_6}$ is given by
\[
H^n(\cZ_{P^2_6}) = \begin{cases}
\F_2 &  n = 0, 8, 9  \\
(\F_2)^{10} & n = 5 \\
(\F_2)^{15} & n = 6 \\ 
(\F_2)^{6} & n = 7 \\ 
0 & \mbox{ otherwise. } 
\end{cases}
\]
Generators in degrees $8$ and $9$ are suspensions of generators in degrees $1$ and $2$ of $H^*(P^2_6)$ and there is a non-trivial $Sq^1$ action taking the degree $1$ generator to the degree $2$ generator. Hence, there is a non-trivial $Sq^1$ action on $\cZ_{P^2_6}$ taking the degree $8$ generator to the degree $9$ generator. All other Steenrod operations act trivially on $H^*(\cZ_{P^2_6})$ since $S^1$ has trivial $Sq^n$ action for $n\geq 1$.

Note that the description of the cohomology of $\cZ_{P^2_6}$ comes from splitting~\eqref{Unstable_splitting} of topological spaces. Given a cochain in $\cZ_{P^2_6}$, one can split it into cochains coming from individual full subcomplexes and vice versa. Thus, combinatorially, knowing how $Sq^1$ acts on these full subcomplexes gives us complete information on how $Sq^1$ acts combinatorially on $\cZ_{P^2_6}$. There are many choices for a simplicial cochain representative of the cohomology class generating $H^1(P^2_6)$. For example, 
 \[x = [1, 4]^* + [1, 6]^* + [2, 5]^* + [2, 6]^* + [4, 5]^*\] %this is x_5
 generates $H^1(P^2_6)$ as does
 \[x^{\prime}= [1, 2]^* + [1, 3]^* + [1, 4]^* + [1, 5]^* + [2, 4]^* + [2, 5]^* + [2, 6]^* + [3, 4]^* + [4, 6]^*.\] % this is x_17
 Here by $[a,b]^{\ast}$, we mean the dual of the corresponding chain $[a,b]$, where
 $[a,b]$ is a $1$-simplex. 
 
Since these classes are of degree one, the operation $Sq^1$ is given by the cup square which at the combinatorial level is concatenation of cochains, that is, $x^2$ has a summand $[a,b,c]^{\ast}$ if and only if the cochain $x$ has summands $[a,b]^{\ast}$ and  $[b,c]^{\ast}$ and $[a,b,c]$ is a simplex in $P^2_6$. Thus, we have
 \[\begin{split}Sq^1(x) &= [1, 4, 5]^*\\
 Sq^1(x^{\prime}) &= [1, 2, 6]^* + [1, 4, 6]^* + [3, 4, 6]^{*}.\end{split}\]
Each of these is a generator of $H^2(P^2_6),$ and so non-trivial, and these expressions give cochain-level representations of 
the non-trivial operation $Sq^1$ on $P^2_6$, which is reflected in the moment-angle complex as the non-trivial operation
\[Sq^1 \colon H^8(\cZ_{P^2_6}) \to H^9(\cZ_{P^2_6}).\]

\subsection{Some results about $Sq^1$ on moment-angle complexes}
The splitting~\eqref{cohlgy_splitting} and Example~\eqref{sec:an_example} allow us to deduce some general results for moment-angle complexes. Recall that the dimension of a simplicial complex is the
maximal dimension of its simplices.

\begin{cor}
    If $\dim(K)=d$, then $Sq^n=0$ in $H^*(\ZK)$ for $n> \lfloor{d/2}\rfloor$. 
    
    In particular, if $K$ is a graph ($d=1$), then there are no non-trivial Steenrod operations on $H^*(\ZK)$.
\end{cor}

\begin{proof}
Recall the instability condition for Steenrod operations acting on the cohomology of a space: a class in degree $i$ only supports non-zero Steenrod operations $Sq^j$ for $0\leq j\leq i$. Thus for any space $X$ whose $\F_2$-cohomology is concentrated in degrees at most $d$, we have $Sq^n=0$ in $H^*(X)$ for $n> \lfloor{d/2}\rfloor$. But if $\dim(K)=d$, this holds for the cohomology of $K$ and all of its full subcomplexes. And then it follows from splitting~\eqref{cohlgy_splitting} 
that the same holds in the cohomology of $\mathcal{Z}_K$.
\end{proof}

\begin{rem}
    Analogous results hold for real moment-angle complexes due to a similar homotopy decomposition after suspension, 
    see~\cite[Corollary 8.3.6(a)]{buchstaber2014toric}. And the same considerations apply more generally for polyhedral products $(\underline{X}, \underline{A})^{K}$ where each $X_i$ is 
    contractible, again because of the homotopy decomposition of the suspension~\cite[Theorem 8.3.5]{buchstaber2014toric}. 
\end{rem}
\begin{cor}
    If $K$ contains $P^2_6$ (the $6$-vertex model of $\mathbb{R}P^2$ in Figure~\ref{rptri}) as a full subcomplex, then there is a non-trivial
$Sq^1$ in $H^*(\mathcal{Z}_K)$.
\end{cor}

\begin{proof}
 This follows from the splitting and the fact that we have a non-trivial $Sq^1$ in $H^*(\mathbb{R}P^2)$.   
\end{proof}

\begin{cor}
    $Sq^1$ acts trivially on classes of dimension $\leq 7$ in an arbitrary $\mathcal Z_K$.
    This bound is optimal.
\end{cor}
\begin{proof}
    Rewriting the isomorphism~\eqref{cohlgy_splitting}
     by incorporating suspensions as shifts of the graded cohomology we get
    \[\widetilde{H}^k(\ZK)
\cong\bigoplus_{J\notin K} \widetilde{H}^{k-|J|-1}(\KJ).\]
Now a class of dimension at most $7$ in an arbitrary $\mathcal Z_K$ is a sum  of classes of dimension $7-|J|-1$  in reduced cohomology of full-subcomplexes, where they are zero when $|J|>6$. For $|J| = 6$, they can be non-zero only when $\KJ$ is disconnected, i.e., a disjoint union of simplicial complexes on fewer vertices. For such a $\KJ$, the cohomology is a direct sum of cohomology of its connected component simplicial complexes, hence a non-zero class will further split into classes coming from simplicial complexes on $5$ or fewer vertices (and the Steenrod square respects this). Hence, every class of dimension $\leq 7$ of an arbitrary moment-angle complex splits into classes in the reduced cohomology of full subcomplexes on $5$ or fewer vertices. But, simplicial complexes on $5$ or fewer vertices have no non-trivial $Sq^1$. This was verified by enumerating all such simplicial complexes and computing their cohomologies using packages for finite simplicial complexes and for cohomology with a basis in Sage \cite{sagemath}, implemented by John H. Palmieri and Travis Scrimshaw, which proves the first statement.

The example of $P^2_6$ shows that this bound is optimal.
\end{proof}

Note that the Sage computation also implies that $P^2_6$ is the minimal triangulation of $\mathbb{R}P^2$.

\begin{rem}
    Define the natural number $m_i$ to be the minimal $m$ for which there exists a simplicial complex $K$ on $m$ vertices admitting a non-trivial $Sq^i$ in mod $2$ cohomology. Then the same argument shows that
    $Sq^i$ acts trivially on classes of dimension $\leq m_i+i$ in an arbitrary $\mathcal Z_K$. The above is the case $i=1$ with $m_1=6$. We are not aware of other known values of $m_i$.
\end{rem}

\setcounter{section}{4} %%%%%%%%%%%%%%%%%%%%%%%%%%%%%%%%%%%%%%%%%%%%%%%%%%%%%%%%%%%
\section{Steenrod operations on cochains of simplicial complexes}\label{sec:combinatorial_formula}
%%%%%%%%%%%%%%%%%%%%%%%%%%%%%%%%%%%%%%%%%%%%%%%%%%%%%%%%%%%

There are several ways of describing Steenrod operations on simplicial cochains combinatorially.  See, for example, 
\cite[Definition 2, Theorem 10]{MEDINAMARDONES2023101921} 
or~\cite[Corollary 3.2]{GONZALEZDIAZ199989}. We give another description which is essentially a reinterpretation of formulas from \cite{MEDINAMARDONES2023101921}. Our interest in this stems from the idea of seeing how such formulas may 
be carried over from a simplicial complex to constructions such as polyhedral products and joins.

Let $K$ be a simplicial complex and $j\in \mathbb N$.
We write $K_j$ for the set of $j$-simplices in $K$, and we write $C_j (K)$, for the simplicial chains  on $K$ with coefficients in $\mathbb{F}_2$.  That is, $C_j(K)$ is the free $\mathbb{F}_2$ vector space generated by the elements of  $K_j$. Similarly, we write $C^j(K)$ for the simplicial cochains
$\mathrm{Hom} (C_j(K), \mathbb{F}_2)$, and we write 
$ \langle -, -\rangle\colon C^j(K)\otimes C_j(K)\to \mathbb{F}_2$
for the pairing of cochains and chains.

Our simplicial complex comes with natural maps called face maps $\partial_i: K_{j} \to K_{j-1}$ for $ 0 \leq i \leq j$ 
which send a $j$-simplex to its $i^{th}$ face. These maps extend linearly to maps $\partial_i \colon C_{j}(K)\rightarrow C_{j-1}(K)$ on the chain level.

\subsection{Combinatorial formula for cochain-level $Sq^1$}

We present first a particularly accessible  description of $Sq^1$ and then show that this does in fact agree with the formula in \cite{MEDINAMARDONES2023101921} before  proceeding to higher $Sq^n$.

For $c \in C^j (K), x \in C_{j+1}(K)$,
\begin{equation}
\label{Sq1-chains}
Sq^1(c)(x) =\sum_{\begin{matrix}0\le u < v \le j+1 \\ u+v \text{ even}\end{matrix}} \mu(\langle c, \partial_u x\rangle  \otimes \langle c, \partial_v x\rangle)
\end{equation}
where  $\mu \colon \mathbb F_2\otimes\mathbb F_2\rightarrow \mathbb F_2$ is the multiplication in $\mathbb F_2$.

Notice that the set $U$ from \cite[Definition 7]{MEDINAMARDONES2023101921} is exactly $U = \{u,v\}$, so the formula (\ref{Sq1-chains}) follows directly from \cite[Definition 2, Theorem~10]{MEDINAMARDONES2023101921}.

While this gives a formula for
$Sq^1$ via the pairing of cochains and chains, 
we can now give an explicit formula internal to cochains by using {\it pasting functions} $f_{u,v}$.

Let $x = [x_0, x_1, \dots, x_j]^*$ and $y = [y_0, y_1, \dots, y_j]^*$ represent generators of $C^j(K)$ with $S = \text{support}(x) \cup \text{support}(y)$ (i.e., $S = \{ x_0, \dots, x_j, y_0, \dots, y_j\}$). If $|S| = j+2$, define $(z_0, \dots, z_{j+1})$ as the elements of $S$ listed in their natural order.

On generators define $f_{u,v} \colon C^j(K) \times C^j(K) \to C^{j+1}(K)$ as
$$f_{u,v}(x,y) = \begin{cases} 
[z_0, \dots, z_{j+1}]^* & \text{if $|S| =  j+2$, $S\setminus\text{support}(x) = \{z_v\}$}\\
&\quad\text{and 
$S\setminus \text{support}(y) = \{z_u\}$}\\
0 &\text{otherwise}
\end{cases}
$$
 and extend $f_{u,v}$ to $C^j(K)\times C^j(K)$ bilinearly.

Finally, for $c\in C^j(K)$,
$$Sq^1(c) = \sum_{\begin{matrix}0\le u < v \le j+1 \\ u+v \text{ even}\end{matrix}} f_{u,v}(c,c).$$ 

\subsection{Combinatorial formula for cochain-level $Sq^n$}
For $n>1$, \cite[Definition 2, Theorem~10]{MEDINAMARDONES2023101921} gives us a similar combinatorial formula for $Sq^n$.

Let $A \subset B \subset \mathbb N$ be ordered finite  sets and let $\mathrm{pos}_B \colon A\rightarrow \mathbb N$ be the function that gives the position of element $a\in A$ in set $B$ (for example if $B = (3,6,7,9)$, then $\mathrm{pos}_B(6) = 2$).

Let $j\in \mathbb N$. For $c \in C^j (K), x \in C_{j+n}(K)$,
\begin{equation}
\label{Sqn-chains}
Sq^n(c)(x) =\sum \mu(\langle c, \partial_{u_1}\partial_{u_2}\cdots\partial_{u_n} x\rangle  \otimes \langle c, \partial_{v_1}\partial_{v_2}\cdots\partial_{v_n} x\rangle)\end{equation}
where $\partial_i \colon C_{j+k}(K)\rightarrow C_{j+k-1}(K)$, $1 \le k \le n$, $0 \le i \le j+k$ is the face operator and $\mu\colon\mathbb F_2\otimes\mathbb F_2\rightarrow \mathbb F_2$ is the multiplication in $\mathbb F_2$. The sum in (\ref{Sqn-chains}) goes over all ordered subsets $U=(u_1,\ldots,u_n)$ and $V=(v_1,\ldots, v_n)$ of $\{0,\ldots,j+n\}$ (i.e., $u_1<u_2<\cdots<u_n$ and $v_1<v_2<\cdots<v_n$) such that $U\cap V = \emptyset$, and for all $i\in \{1,\ldots, n\}$
\[\mathrm{pos}_{U\cup V}(u_i)\equiv_2 u_i,\;  \mathrm{pos}_{U\cup V}(v_i)\not\equiv_2 v_i.  \]

Notice sets $U$ and $V$ are exactly sets $U^0$ and $U^1$ from \cite[Definition 7]{MEDINAMARDONES2023101921}. The only difference is that sets $U$ and $V$ have $n$ elements each and $U^0$ and $U^1$ have $2n$ elements in total, but the end result for $Sq^n(c)(x)$ is the same since all summands in (\ref{Sqn-chains}) for $U^0$ and $U^1$ not having the same number of elements will vanish. Once again, the formula (\ref{Sqn-chains}) follows from \cite[Definition 2, Theorem~10]{MEDINAMARDONES2023101921}.

We now give an explicit formula on cochains by using {\it pasting functions} $f_{U,V}$.

Let $x = [x_0, x_1, \dots, x_j]^*$ and $y = [y_0, y_1, \dots, y_j]^*$ represent generators of $C^j(K)$ with $S = \text{support}(x) \cup \text{support}(y)$ (i.e., $S = \{ x_0, \dots, x_j, y_0, \dots, y_j\}$). If $|S| = j+n+1$, define $(z_0, \dots, z_{j+n})$ as the elements of $S$ listed in their natural order.  Let $U,V\subset \{0,\ldots,j+n\}$.

On generators define $f_{U,V} \colon C^j(K) \times C^j(K) \to C^{j+n+1}(K)$ as
$$f_{U,V}(x,y) = \begin{cases} 
[z_0, \dots, z_{j+n}]^* & \text{if $|S| =  j+n+1$, $S\setminus\text{support}(x) = V$}\\
&\quad\text{and 
$S\setminus \text{support}(y) = U$}\\
0 &\text{otherwise}
\end{cases}
$$
 and extend $f_{U,V}$ to $C^j(K)\times C^j(K)$ bilinearly.

Finally, 
\begin{equation}
\label{Sqn-cochains}
    Sq^n(c) = \sum f_{U,V}(c,c),
    \end{equation}
where the sum in (\ref{Sqn-cochains}) goes over all ordered subsets $U=(u_1,\ldots,u_n)$ and $V=(v_1,\ldots, v_n)$ of $\{0,\ldots,j+n\}$ such that $U\cap V = \emptyset$, and for all $i\in \{1,\ldots, n\}$
\[\mathrm{pos}_{U\cup V}(u_i)\equiv_2 u_i, \;  \mathrm{pos}_{U\cup V}(v_i)\not\equiv_2 v_i.  \]
\smallskip

In the next section we will initiate the study of how these combinatorial formulas
for Steenrod operations of a simplicial complex $K$ may be used to 
 investigate the effect of combinatorial changes in $K$ on the Steenrod action. In particular this leads us to consider the construction known as polyhedral join.

\setcounter{section}{5} %%%%%%%%%%%%%%%%%%%%%%%%%%%%%%%%%%%%%%%%%%%%%%%%%%%%%%%%%%%
\section{Polyhedral joins}
%%%%%%%%%%%%%%%%%%%%%%%%%%%%%%%%%%%%%%%%%%%%%%%%%%%%%%%%%%%

In this section we introduce polyhedral joins and investigate Steenrod operations on their cohomology. We begin with special cases where the homotopy type is understood as that of a join. Then we consider substitution complexes for which we establish a homotopy decomposition. This is then used to establish how non-trivial Steenrod operations on the building complexes are propagated to the substitution complex.

\begin{defn}
{Given a simplicial complex $K$ on $m$ vertices, and $m$ topological pairs $(\underline X, \underline A)=((X_1, A_1),\ldots, (X_m, A_m))$, the \emph{polyhedral join} $(\underline X, \underline A)^{*K}$ is the subspace of the join
${\mathlarger{\mathlarger{\mathlarger{\mathlarger{\mathlarger{\ast}}}}}}_{i =1}^m X_i$ given by
\[
(\underline X, \underline A)^{*K}=\bigcup_{\sigma\in K} (\underline X, \underline A)^\sigma\]
where \[(\underline X, \underline A)^\sigma=\overset{m}{\underset{i=1}{\mathlarger{\mathlarger{\mathlarger{\mathlarger{\mathlarger{\ast}}}}}}} Y_i\quad {\rm \ and \ }\quad Y_i = \begin{cases} X_i &{\rm \ if \ } i \in \sigma \\ A_i & {\rm \ if \ } i \notin \sigma.\end{cases}\]}
\end{defn}

Here we see that if $K$ is an $(m-1)$-simplex then $(\underline X, \underline A)^{*K}= X_1 * \dots * X_m$, and if $K$ is a disjoint union of $m$ points then $(\underline X, \underline\emptyset)^{*K} = X_1 \sqcup \dots \sqcup X_m.$

While the polyhedral join is defined for topological pairs $(\underline{X}, \underline{A})$, in this paper we are mostly concerned with this construction for pairs $(\underline{K}, \underline{L})$ of simplicial complexes.  In this case, $(\underline{K}, \underline{L})^{*K}$ is again a simplicial complex, not just a topological space.

\subsection{Special cases}
The following table records some other important cases of the polyhedral join.
\begin{center}
\renewcommand{\arraystretch}{1.3}
\begin{tabular}{l|c|l}
 $(\underline{K}, \underline{L})$&polyhedral join $(\underline{K},\underline{L})^{*K}$&name
 \\
 &(notation)&\\
 \hline\hline
$(\underline{\Delta^{1}}, \underline{\partial\Delta^{1}})$ &$D(K)$&simplicial double\\
\hline
$(\underline{\Delta^{l}}, \underline{\partial\Delta^{l}})$
&$K(\Delta^{l_1}, \dots, \Delta^{l_m})$&simplicial wedge\\
\hline
$(\underline{\Delta^{l}}, \underline{L})$
&$K(L_1, \dots, L_m)$&simplicial composition
\\
\end{tabular}
\end{center}
\smallskip

Each line of the table above is successively more general
and the homotopy type is known to be a join~\cite[Proposition 6.1]{A13}
\begin{equation} \label{heq-join}
K(L_1, \dots, L_m)\simeq K*L_1*\dots *L_m.
\end{equation}
Since a join is homotopy equivalent
to the suspension of the smash product, its reduced cohomology is given by the shifted tensor product of the cohomologies. This is also recovered as~\cite[Example 3.11]{Z18}.

\begin{prop}
\label{prop:splitting_join}
There is an isomorphism of $\mathcal{A}$-modules
\[
\widetilde{H}^*(K(L_1, \dots, L_m))
\cong \Sigma^m \widetilde{H}^*(K)\otimes 
\widetilde{H}^*(L_1)\otimes \dots \otimes \widetilde{H}^*(L_m).
\]  
\end{prop}

\begin{proof}
This follows directly from the homotopy equivalence~\eqref{heq-join}, using that $X*Y\simeq \Sigma X\wedge Y$, that $\widetilde{H}^*(X\wedge Y)\cong \widetilde{H}^*(X)\otimes \widetilde{H}^*(Y)$ and the stability of Steenrod operations.
\end{proof}

\subsection{Substitution complexes}

    Let $K$ be a simplicial complex on $[m]$ and $K_i$, $i\in[m]$, simplicial complexes on $[n_i]$. The {\it substitution complex} is defined in \cite{AbramyanPanov} as the polyhedral join
    \[K\langle K_1,\ldots,K_m\rangle = (\underline{K},\underline{\emptyset})^{*K}.\]

We can label vertices of $K\langle K_1, \ldots, K_{m}\rangle$ by $v_{i,j}$, $i\in [m]$ and $j\in [n_i]$ to indicate that they arise in $K\langle K_1, \ldots, K_{m}\rangle$ as a result of the substitution of the complex $K_i$ into the vertex $v_i \in K$. 

In what follows, let $\bullet$ denote a simplicial complex consisting of a single vertex. We will also use the following notation.
\begin{equation}   \label{def:K^i}
\begin{split}
        &K^{0}=K,\\
        &K^1=K\langle K_1, \bullet, \ldots, \bullet \rangle,\\
        &K^2=K^1\langle \underbrace{\bullet, \ldots, \bullet,}_{n_1} K_2, \bullet, \ldots, \bullet \rangle\\
        &\;\;\;\vdots \\
        &K^m=K^{m-1}\langle\underbrace{\bullet, \ldots, \bullet,}_{\sum_{i=1}^{m-1} n_i} K_m \rangle.
    \end{split}
    \end{equation}

Observe that, by the definition of the substitution complex, the complex $K$ is a full subcomplex of $K\langle K_1, \ldots, K_m\rangle$ on vertices $v_{1,k_1},\ldots,v_{m,k_m}$ for any $k_{j}\in [n_j]$. In particular, we can choose $k_j=1$ for all $j$, then, when we talk about the vertex $v_i$, $i\in[m]$ in $K\langle K_1,\ldots,K_m\rangle$, we identify $v_i$ with the vertex $v_{i,1}$. 

The main result of this section is the homotopy decomposition of substitution complexes described in 
Theorem~\ref{thm:substitution_cpx_splitting} below.
For the statement of the theorem and for further discussion, we recall the definition of a link.    Let $K$ be a simplicial complex and $v\in K$ a vertex. The {\it link of $v$ in $K$} is the subcomplex
    \[\link_K(v) = \{\tau\in K\mid v\notin \tau, \tau\cup\{v\}\in K\}.\]
    We omit the subscript $K$ and write simply $\link(v)$ when it is obvious in which simplicial complex the link is taken.

\begin{thm}\label{thm:substitution_cpx_splitting}
    Let $K, K_1,\ldots, K_m$ be simplicial complexes on $[m]$, $[n_1],\ldots, [n_m]$, respectively. Assume that $K$ is connected. Then 
    \[
K\langle K_1, \ldots, K_m\rangle\simeq K\vee \link_{K^0}(v_1)*K_1\vee\ldots\vee \link_{K^{m-1}}(v_m)*K_m.
\]
\end{thm}

To prove the theorem, we first prove several technical statements.

\begin{prop}\label{prop:sequential_substitution}
    The substitution complex $K\langle K_1, \ldots, K_m \rangle$ is isomorphic to the substitution complex
    \begin{equation*}
        (K\langle K_1, \ldots, K_{m-1}, \bullet \rangle) \langle \underbrace{\bullet, \ldots, \bullet }_{n_1}, \ldots,  \underbrace{\bullet, \ldots, \bullet }_{n_{m-1}}, K_m    \rangle .
    \end{equation*}
    It follows that a substitution complex can be constructed sequentially.
\end{prop}
\begin{proof}
    Let $\sigma\in K\langle K_1,\ldots,K_{m-1},\bullet\rangle$ be a simplex not containing the vertex $v_m$. Observe that we can represent $\sigma$ as $\sigma = \tau_{i_1}*\cdots*\tau_{i_k}$, where $\tau_{i_j} \in K_{i_j}$, for $\{v_{i_1},\ldots,v_{i_k}\}\in K$ and $m\notin \{i_1,\ldots,i_k\}$. Next note that simplices of 
    \begin{equation}
        \label{sequen}
     (K\langle K_1, \ldots, K_{m-1}, \bullet \rangle) \langle \underbrace{\bullet, \ldots, \bullet }_{n_1}, \ldots,  \underbrace{\bullet, \ldots, \bullet }_{n_{m-1}}, K_m   \rangle
         \end{equation}
     are all possible $\sigma\in K\langle K_1,\ldots,K_{m-1},\bullet\rangle$ not containing the vertex $v_m$ and additionally all simplices $\sigma*\tau_m$ for $\tau_m\in K_m$, whenever $\{v_{i_1},\ldots,v_{i_k},v_m\}\in K$. In other words, all simplices in~\eqref{sequen} are of the form $\tau_{i_1}*\cdots*\tau_{i_k}$ for $\{v_{i_1},\ldots,v_{i_k}\}\in K$ (here we allow $m$ to be in $\{i_1,\ldots,i_k\}$) and $\tau_{i_j}\in K_j$. These are exactly the simplices of $K\langle K_1,\ldots,K_m\rangle$, so the claim holds.
\end{proof}

   The {\it star of $v$ in $K$} is the subcomplex 
    \[\star_K(v) = \{\tau\in K\mid \tau\cup\{v\}\in K\}.\]

    We omit the subscript $K$ and write $\star(v)$ when $K$ is understood from the context.

\begin{prop}\label{prop:gluing_subst_cpx}
    Let $K$ and $K_i$ be simplicial complexes and $v_i \in K$ a vertex. Then 
    \[
    K\langle \bullet, \ldots, \bullet, K_i,\bullet \ldots, \bullet\rangle \cong K \bigsqcup_{\star_K(v_i)} \star_K(v_i)\langle K_i, \underbrace{\bullet, \ldots, \bullet}_{s_i-1} \rangle
    \]
    where the isomorphism is of simplicial complexes, $s_i$ is the number of vertices in $\star(v_i)$ and, for the ease of notation, we assume that $v_i$ is the first vertex of $\star(v_i)$.
\end{prop}

\begin{proof}
    When substituting $K_i$ in $K$, we are placing $K_i$ instead of the vertex $v_i$ and joining it with everything $v_i$ was connected to, i.e., we are joining it with $\link_K(v_i)$. This means we have an isomorphism
    \[
        K\langle \bullet, \ldots, \bullet, K_i,\bullet \ldots, \bullet\rangle \cong K|_{[m]\backslash \{v_i\}} \bigsqcup_{\link(v_i)} \link(v_i)*K_i.\]
    Similarly $\link(v_i)*K_i$ can be seen as $\star(v_i)\langle K_i,\underbrace{\bullet, \ldots, \bullet}_{s_i-1} \rangle$, where $s_i$ is the number of vertices of $\star(v_i)$. Thus 
     \[\begin{split}K\langle \bullet, \ldots, \bullet, K_i,\bullet \ldots, \bullet\rangle &\cong K|_{[m]\backslash \{v_i\}} \bigsqcup_{\link(v_i)} \star(v_i)\langle K_i,\underbrace{\bullet, \ldots, \bullet}_{s_i-1} \rangle\\
     &\cong K \bigsqcup_{\star(v_i)} \star(v_i)\langle K_i, \underbrace{\bullet, \ldots, \bullet}_{s_i-1} \rangle.
     \end{split} 
     \] 
\end{proof}

\begin{lem}\label{lem:wedge_of_cw_complexes}
Let $K, L$ and $M$ be connected $CW$-complexes and $L\simeq M$. Then $(K,v)\vee (L,v') \simeq (K,v)\vee (M,v'')$ for any choice of base points $v\in K, v'\in L$ and $v''\in M$.    
\end{lem}
\iffalse
\textcolor{brown}{Evgeniya: I updated the proof. If this seems too long, i think it can be shortened to the first paragraph plus: 
"The lemma follows from the application of Proposition C.3.5 in \cite{buchstaber2014toric} to push-outs $(K,v)\vee (L,v')$ and $(K,v)\vee (M,v'')$. "}
\fi
\begin{proof}
    Let $\varphi\colon L\rightarrow M$ and $\psi\colon M\rightarrow L$ be such that $\varphi\circ\psi\simeq id_M$ and $\psi\circ\varphi\simeq id_L$. Since the spaces are path-connected and well-pointed for any choice of base point, we can assume that $\varphi(v')=v''$ and $\psi(v'')=v'$.

    The space $(K,v)\vee (L,v')$ is a push-out of
    \[
    K\leftarrow \{v\} \rightarrow L
    \]
    and similarly, the space $(K,v)\vee (M,v'')$ is a push-out of 
     \[
    K\leftarrow \{v\} \rightarrow M.
    \]
    Since $CW$-complexes are cofibrant objects and the inclusions of $v$ in the above diagrams are cofibrations, we apply ~\cite[Proposition C.3.5]{buchstaber2014toric} to conclude that these push-out objects are weakly (and, thus for $CW$-complexes also strongly) homotopy equivalent to their homotopy colimits. Since we also have homotopy equivalences between the respective objects of the diagrams, we conclude that the push-outs are homotopy equivalent.
\end{proof}

\begin{lem}\label{lem:homotopy_type_one_substituion}
    Let $K$ and $K_i$ be simplicial complexes and $v_i \in K$ a vertex. Assume that $K$ is connected. Then 
    \[
    K\langle \bullet, \ldots, \bullet, K_i,\bullet \ldots, \bullet \rangle \simeq (K,v_i) \vee (\star(v_i)\langle K_i,\bullet,\ldots,\bullet\rangle,v_i),
    \]
    
    where we assume that $v_i$ is the first vertex of $\star(v_i)$.
\end{lem}
\begin{proof}
Let $S =  K\langle \bullet, \ldots, \bullet, K_i,\bullet \ldots, \bullet \rangle$. Since $(S,\star(v_i))$ has the homotopy extension property and $\star(v_i)$ is contractible, we have a homotopy equivalence of topological spaces
\[S\simeq S/\star(v_i).\]
From Proposition \ref{prop:gluing_subst_cpx}, we now get
\[\begin{split}
   S&\simeq (K/\star(v_i)) \vee (\star(v_i)\langle K_i,\bullet,\ldots,\bullet\rangle / \star(v_i))\\
   & \simeq K\vee \star(v_i)\langle K_i,\bullet,\ldots,\bullet\rangle
   \end{split}\]
where the last homotopy equivalence follows from Lemma \ref{lem:wedge_of_cw_complexes}. 

Since both $K$ and $\star(v_i)\langle K_i,\bullet,\ldots,\bullet\rangle$ are path connected, we can choose $v_i$ to be the base point for both of them, i.e.,
\begin{equation}\label{hom_eq_S}
    S \simeq (K,v_i) \vee (\star(v_i)\langle K_i,\bullet,\ldots,\bullet\rangle,v_i). \qedhere
\end{equation}
\end{proof}

We are now ready to prove Theorem~\ref{thm:substitution_cpx_splitting}.

\begin{proof}[Proof of Theorem~\ref{thm:substitution_cpx_splitting}]
    Recall the definition (\ref{def:K^i}) of simplicial complexes $K^0,K^1,\ldots, K^m$.    By Proposition~\ref{prop:sequential_substitution}, we have that $K^m\cong K\langle K_1,\ldots,K_m \rangle$ and, by the definition of the substitution complex, we have that the complex $K^i$ is a full subcomplex of $K^{i+1}$ on vertices
    \begin{align*}
    &v_{1,1}, \ldots, v_{1,n_1},\\
    &\vdots \\
    &v_{i,1}, \ldots, v_{i,n_i},\\
    &v_{i+1,1},\\
    &\vdots\\
    &v_{m,1} .
    \end{align*}
By Lemma~\ref{lem:homotopy_type_one_substituion}, we have a homotopy equivalence
$$K^{i+1}\simeq (K^{i}, v_{i+1, 1})\vee (\star_{K^i}(v_{i+1, 1})\langle K_{i+1}, \bullet, \ldots, \bullet \rangle, v_{i+1, 1}).$$
To simplify the notation, let $X_{i+1}=\star_{K_i}(v_{i+1, 1})\langle K_{i+1}, \bullet, \ldots, \bullet \rangle$. Notice that $X_{i+1} = \link_{K^i}(v_{i+1,1})*K_{i+1}$.

Thus, to get to the conclusion of the theorem, we observe the sequence of homotopy equivalences
\begin{align*}
    K\langle K_1, \ldots, K_m\rangle \cong K^m &\simeq (K^{m-1}, v_{m, 1})\vee (X_{m}, v_{m, 1})\\
    &\simeq \Big((K^{m-2}, v_{m-1, 1})\vee(X_{m-1}, v_{m-1, 1}), v_{m, 1} \Big) \vee (X_{m}, v_{m, 1})\\
    &\vdots\\
    &\simeq \Big(\cdots\Big(\Big((K^{0}, v_{1, 1})\vee(X_{1}, v_{1, 1}), v_{2, 1}\Big)\vee(X_{2}, v_{2, 1})\Big) \ldots \Big) \vee (X_{m}, v_{m, 1})
\end{align*}
which we can obtain by Lemma~\ref{lem:homotopy_type_one_substituion}.

The claim follows since $K$ is connected, by an application of Lemma~\ref{lem:wedge_of_cw_complexes} to change the base point to $v_{1,1}$.
\end{proof}

\subsection{Propagation of Steenrod operations on substitution complexes}
In this section we consider situations where we can deduce the existence of a non-trivial 
Steenrod operation on a substitution complex from the existence of such an operation for one of the simplicial complexes used to build the substitution complex. We refer to this scenario as \emph{propagation} of Steenrod operations. 

Here we establish some results on how non-trivial Steenrod operations on $K$ and the $K_i$s propagate to the substitution complex $K\langle K_1,\ldots, K_m\rangle$ given the splitting from Theorem~\ref{thm:substitution_cpx_splitting}. Further we study explicitly how these propagations manifest combinatorially with examples. 

\begin{prop}\label{prop:non-trivial-Sqn-substitution}
    If $K$ has a non-trivial $Sq^n$, then the same is true for the substitution complex $K\langle K_1,\ldots, K_m\rangle$.
\end{prop}
 \begin{proof}
     This is immediate from Theorem~\ref{thm:substitution_cpx_splitting}. The simplicial complex $K$ appears as a wedge summand in the splitting of the substitution complex and thus $\widetilde{H}^{\ast}(K)$ is a summand in the cohomology of the substitution complex.
 \end{proof}

By similar reasoning, we get the following propagation result.

\begin{prop}
    If $K_i$ has a non-trivial $Sq^n$ and $\link_{K}(v_i)$ is not acyclic, then the substitution complex $K\langle K_1,\ldots, K_m\rangle$ has a non-trivial $Sq^n$ propagating from $K_i$.
\end{prop}
\begin{proof}
    Reorder the vertices of the complex $K$, so that $i$ is now the first vertex (and in the substitution complex $K_i$ is now $K_1$). 
    By Theorem~\ref{thm:substitution_cpx_splitting} for this reordered substitution complex, up to suspension, $\widetilde{H}^{\ast}(\link_{K}(v_i)) \otimes \widetilde{H}^{\ast}(K_i)$ is a summand in the cohomology of the substitution complex. Using non-acyclicity of $\link_{K}(v_i)$ and stability of Steenrod operations, we have the result.
\end{proof}

In general, one cannot propagate $Sq^n$ from $K_i$ to the substitution complex $K\langle K_1,\ldots, K_m\rangle$. 
\begin{eg}
Let $K=\Delta^1$ and let $K_1=P^2_6$. Then $K\langle K_1, \bullet\rangle\cong \Cone P^2_6\simeq *$ does not inherit a non-trivial $Sq^1$ from $P^2_6$. More generally, for $K=\Delta^1$, the substitution complex $K\langle K_1, \bullet\rangle=\Cone K_1\simeq *$. 
A whole family of examples of this type can be derived by requiring, for example, that $\mathrm{link}(v_i)$ is contractible. 
\end{eg}

\begin{rem}
    Each $K_i$ is a full subcomplex of the substitution complex 
    $K\langle K_1,\ldots, K_m\rangle$, but this does not ensure a non-trivial propagation of the Steenrod operations (unlike in the case of moment-angle complexes). We return to this below.
\end{rem}

\subsubsection{Combinatorial study of propagation in substitution complexes}
By Proposition \ref{prop:non-trivial-Sqn-substitution},
a non-trivial Steenrod action in $K$ propagates to the substitution complex $K\langle K_1,\ldots,K_m\rangle$. If we know a cochain representative of a class in $H^*(K)$ on which the Steenrod action is non-trivial, the question arises as to how to obtain a representative of a class in the substitution complex with non-trivial Steenrod action. One can then use the formulas from Section~\ref{sec:combinatorial_formula} to calculate the non-trivial Steenrod operation on this class.

To begin, let us focus on a simple example.
\begin{eg}
Let $K = P^2_6$, $K_1 = S^0$ and $K_i = \bullet$, $i\in\{2,3,4,5,6\}$.

Here, a non-trivial Steenrod square is $Sq^1$ from the first to the second cohomology group of $K$. To address our question, we want to know which conditions  ensure that a cochain is a cocycle and not a coboundary, i.e., that it represents a non-trivial cohomology class in the substitution complex  $P^2_6\langle S^0,\bullet, \bullet,\bullet,\bullet,\bullet\rangle$. 

For a simplicial complex $K$, in order for a one-dimensional cochain to be a cocycle, it must contain an odd length cycle of duals of edges. And in order for it to not be a coboundary, for each 2-simplex in $K$, this cochain needs to contain the duals of an even number of its edges (so either 0 or 2). For higher dimensions, we do not know the corresponding conditions. For $P^2_6$ this works out nicely. First, since there are no 3-simplices, each two-dimensional cochain is a cocycle. Second, since each edge belongs to exactly two 2-simplices, if a two-dimensional cochain is the sum of an odd number of duals of 2-simplices, it is not a coboundary.

Recall from Section \ref{sec:an_example} a cochain representative
\[x = [1, 4]^* + [1, 6]^* + [2, 5]^* + [2, 6]^* + [4, 5]^*\] %this is x_5
of the generator in $H^1(P^2_6)$.

The substitution complex $P^2_6\langle S^0,\bullet, \bullet,\bullet,\bullet,\bullet\rangle$ is obtained by adding a new vertex $7$ and connecting it to everything vertex $1$ was connected to, i.e., we are adding facets $[2,3,7]$, $[2,6,7]$, $[3,5,7]$, $[4,5,7]$ and $[4,6,7]$ to $P^2_6$.

In order to obtain a representative of the generator in the substitution complex, we start from $x$ and then add edges to it  to construct a cocycle that is not a coboundary. It is enough for each edge in $x$ containing the vertex 1 to add an edge where 1 is exchanged for 7. Thus we get
\[y = [1, 4]^* + [1, 6]^* + [2, 5]^* + [2, 6]^* + [4, 5]^* + [4,7]^* + [6,7]^* = x + [4,7]^* + [6,7]^*.\]

Now we can use the combinatorial formula in order to calculate $Sq^1(y)$, which in this case is the cup product, i.e. concatenation of cochains. Thus, $y^2$ has a summand $[a,b,c]^{\ast}$ if and only if the cochain $y$ has summands $[a,b]^{\ast}$ and $ [b,c]^{\ast}$ and $[a,b,c]$ is a simplex in $P^2_6\langle S^0,\bullet, \bullet,\bullet,\bullet,\bullet\rangle$. This gives us
\[Sq^1(y) = y^2 = [1,4,5]^* + [2,6,7]^*.\]
Since Steenrod squares commute with the coboundary operator, $Sq^1(y)$ is a cocycle and it is not difficult to check that it is not a coboundary, so $y$ is indeed a representative of a class with non-trivial Steenrod action.
\end{eg}

Next we consider replacing $S^0$ with $S^1$ in the example above, where $S^1$ is modelled as the simplicial complex $\partial\Delta^2$.

\begin{eg}
Let $K = P^2_6$, $K_1 = \partial\Delta^2$ and $K_i = \bullet$, $i\in\{2,3,4,5,6\}$.

The substitution complex $P^2_6\langle \partial\Delta^2,\bullet, \bullet,\bullet,\bullet,\bullet \rangle$ can be obtained by adding two vertices 7 and 8 and 15 facets $[a,b] * [c,d]$ for each edge $[a,b]$ in $\link_{P^2_6}(1)$ and $[c,d]$ an edge in triangle $[1,7,8]$. If we start again with the same cochain
\[x = [1, 4]^* + [1, 6]^* + [2, 5]^* + [2, 6]^* + [4, 5]^*\] %this is x_5
and for each edge containing vertex 1 we add two new edges where 1 is exchanged for 7 and then 8, we will get a representative of a non-trivial class in the substitution complex. Namely, we take
\[\begin{split}
    y &= [1, 4]^* + [1, 6]^* + [2, 5]^* + [2, 6]^* + [4, 5]^* + [4,7]^* + [4,8]^* + [6,7]^* + [6,8]^*\\ 
    &= x + [4,7]^* + [4,8]^* + [6,7]^* + [6,8]^*.
    \end{split}\]
Now, by the combinatorial formula we have
\[Sq^1(y) =y^2 = [1,4,5]^* + [1,4,7]^* + [1,4,8]^* + [1,6,7]^* + [1,6,8]^* + [2,6,7]^* + [2,6,8]^*.\]

Again, since Steenrod squares commute with the coboundary operator, $Sq^1(y)$ is a cocycle. Notice that all edges except for $[1,7]$, $[1,8]$ and $[7,8]$ belong to an even number of 2-simplices, so if a two-dimensional cochain contains an odd number of duals of 2-simplices whose boundary does not have any of the edges $[1,7]$, $[1,8]$ and $[7,8]$, then this cochain is not a coboundary. This is true for $Sq^1(y)$ (it contains three such 2-simplices $[1,4,5]$, $[2,6,7]$ and $[2,6,8]$) and thus, $Sq^1(y)$ is non-trivial.
\end{eg}

The same procedure can be carried out for more general substitutions $P^2_6\langle K_1,\ldots, K_6\rangle$.

In general, if $K$ is an arbitrary simplicial complex with non-trivial Steenrod action, describing appropriate cochains can be more difficult. The problem lies in determining whether a cochain is a cocycle and not a coboundary. With $P^2_6$ we were only dealing with cochains of dimension 1 and 2 where things are still not too complicated, but in higher dimensions it is far more difficult to give these non-triviality conditions. If we were to have a simple enough way to determine whether a cochain is a representative of a non-trivial class we could carry out a similar procedure to that described above. Start with a non-trivial cochain representative $x$ in $H^*(K)$ which supports a non-trivial $Sq^n$. Next, extend it ``minimally'' by adding what is needed for it to become a non-trivial representative $y$ in the substitution complex. Then calculate $Sq^n(y)$ by using the combinatorial formula (\ref{Sqn-cochains}) and finally, check if $Sq^n(y)$ is non-trivial.

\subsection{Propagation of Steenrod operations on moment-angle complexes} The interaction of moment-angle complexes and polyhedral joins of pairs of simplicial complexes has been of interest in many contexts (see for example~\cite{Vidaurre}, \cite{GSS}).
In this section we consider various ways Steenrod operations on the cohomology of $K, K_1,\ldots, K_m, L_1,\ldots, L_m$ can be propagated to Steenrod operations on the cohomology of moment-angle complexes over the polyhedral joins $(\underline{K}, \underline{L})^{*K}$.

Recall that the Hochster formula determines the cohomology of the moment-angle complex $\mathcal Z_K$ in terms of the cohomology of all full subcomplexes $K_J$ of $K$. The homotopy type of a full subcomplex $K_J$ of $K$ is not in general related to the homotopy type of $K$. For an illustration, one can consider for $K$ a cone over an arbitrary simplicial complex $L$ and thus $K$ is contractible, while $L$, which is a full subcomplex of $K$, can be of an arbitrary homotopy type.  We have a slightly better control over the homotopy type of certain polyhedral joins and their particular full subcomplexes.

Consider a composition complex $(\underline{\Delta}, \underline{L})^{*K}$, where $(\underline{\Delta}, \underline{L})=(\Delta^{n_i-1}, L_i)_{i=1}^m$ . As noted in~(\ref{heq-join}), $(\underline{\Delta}, \underline{L})^{*K}$ is homotopy equivalent to $K*L_1*\ldots*L_m$. It is easy to notice that in general neither $K$ nor the $L_i$s are full subcomplexes of the composition complex $(\underline{\Delta}, \underline{L})^{*K}$. As an example, take for each of $K, L_1$ and $L_2$ two disjoint vertices. Then the composition complex $(\Delta^1,S^0)^{*S^0}$ is the simplicial complex $\partial\Delta^3$, the boundary of a $3$-simplex, and no two disjoint vertices form a full subcomplex of it. 

\begin{lem}
 Let $K$ be a simplicial complex on $[m]$ and $L_1,\ldots, L_m$ simplicial complexes on $n_i$ vertices (allowing ghost vertices). If none of $K$, $L_1,\ldots, L_m$ is contractible, then the non-trivial Steenrod operations on $K$, $L_1,\ldots, L_m$ propagate to non-trivial Steenrod operations on the moment-angle over the composition complex $(\underline{\Delta}, \underline{L})^{*K}$.
\end{lem}
\begin{proof}
By the Hochster formula, the cohomology of the suspended composition complex $(\underline{\Delta}, \underline{L})^{*K}$ appears as a summand of the cohomology of the moment-angle complex. By the Cartan formula and the fact that the composition complex is homotopy equivalent to $K*L_1*\ldots*L_m$, the non-trivial Steenrod operations on $K$, $L_1,\ldots, L_m$ propagate to non-trivial Steenrod operations on the moment-angle complex.
\end{proof}

We end by returning to substitution complexes $K\langle K_1,\ldots,K_m\rangle=(\underline{K}, \underline{\emptyset})^{*K}$.  By Theorem~\ref{thm:substitution_cpx_splitting}, 
\[
K\langle K_1, \ldots, K_m\rangle\simeq K\vee \link_{K^0}(v_1)*K_1\vee\ldots\vee \link_{K^{m-1}}(v_m)*K_m.
\]

In the previous section, we commented that the non-trivial Steenrod operations on $K_i$ do not necessarily propagate to non-trivial Steenrod operations on the substitution complex $(\underline{K}, \underline{\emptyset})^{*K}$, in particular, when $\mathrm{link}(v_i)$ is contractible. Nevertheless, we can propagate the non-trivial Steenrod operations on $K_i$ to those on the moment-angle complex over the substitution complex $(\underline{K}, \underline{\emptyset})^{*K}$.
It is easy to see that all $K$ and $K_1,\ldots, K_m$ appear as full subcomplexes of the substitution complex $(\underline{K}, \underline{\emptyset})^{*K}$.

\begin{prop}
Let $K$ be a simplicial complex on $[m]$ and $K_1,\ldots, K_m$ simplicial complexes on $n_i$ vertices, respectively. Then the non-trivial Steenrod operations on $K$, $K_1,\ldots, K_m$ give rise to non-trivial Steenrod operations on the moment-angle complex over the substitution complex $(\underline{K}, \underline{\emptyset})^{*K}$.
\end{prop}
\begin{proof}
By the Hochster formula, the cohomology of the moment-angle complex over the substitution complex $(\underline{K}, \underline{\emptyset})^{*K}$ decomposes into the direct sum of the cohomology of all (suspended) full subcomplexes of the substitution complex. In the substitution complex itself, $K$, $K_1,\ldots, K_m$ are all full subcomplexes, so all the Steenrod operations on their cohomology will give propagate to non-trivial Steenrod operations on the moment-angle complex.
\end{proof}
\bibliographystyle{abbrv} 
\bibliography{bibliography}
\end{document}